# Analytical solution of general interface problems for reaction-advection-diffusion equations

Andreas Chatziafratis [1,2,3*], Guenbo Hwang [4], Georgia Karali [2,3]

[1] Waseda Research Institute for Science and Engineering, Waseda University, Tokyo, Japan
[2] Department of Mathematics, National and Kapodistrian University of Athens, Greece
[3] Institute of Applied and Computational Mathematics, Foundation for Research and Technology – Hellas, Crete
[4] Department of Mathematics, Daegu University, Gyeongsan, Gyeongbuk, Korea

November 2025

**Abstract.** We obtain novel integral representations, expressed as contour integrals in the complex Fourier plane, for the solution of fully nonhomogeneous interface problems for classical reaction-advection-diffusion (RAD) equations with arbitrary initial data on the line and general interface conditions at the origin. Such problems of applied mathematics emerge in a spectacular variety of physical phenomena in mathematical physics, chemistry, continuum mechanics, biology, pharmacokinetics, as well as mathematical finance. A suitable implementation of Fokas' unified transform method is in force herein, with particular challenges arising due to inclusion in the model of lower-order terms, and the generality of interface problems under consideration. To our knowledge, this careful analysis has not appeared in (nor does it follow easily from) literature published so far. Our explicit formulae directly lend themselves to explorations of solutions' qualitative properties such as regularity, asymptotic behavior and well-posedness. This work has also been of utility in the investigation of nonlinear counterparts as well as in the study of phase-transition phenomena, moving-interface and diffuse-interface problems. Such results for RAD and high-order PDEs will soon appear in a series of separate papers.

## 1. Introduction

Interface problems for partial differential equations (PDEs) arise when distinct physical processes take place in different media, coupled through transmission laws across at their common boundaries. The governing equation is formulated separately in each subdomain, while the interface conditions enforce continuity or balance relations arising from conservation principles. Depending on the modeling framework, these conditions may involve continuity of the solution and its flux, or more general linear transmission relations.

Canonical examples of such models include heat conduction problems in which different media are joined together and the thermal behavior across their boundaries must be properly described. A classical example is the heat conduction problem in a composite wall [1,2]. In this setting, each layer is in perfect thermal contact with constant material properties, and extends infinitely in the directions parallel to the wall. At the interfaces, continuity of temperature and heat flux is imposed, while on the external boundaries one may prescribe Dirichlet, Neumann, or Robin conditions. This framework also extends naturally to one-dimensional rods in thermal contact. From a physical perspective, another relevant case is heat conduction with frictional heating, which arises when moving bodies are in thermal contact [3,4]. In this context, the interface is the contact surface where friction generates localized heat flux, and the resulting heat flow is distributed across the surfaces such that their sum equals the specific power of friction.

Applications of such interface problems are widespread, including pollutant transport in stratified aquifers, heat conduction in composite strips, reaction-diffusion processes in layered biological tissues, chemical transport in heterogeneous catalysts and so on [4-12]. In these settings, discontinuities in the coefficients represent variations in material properties, while the interface conditions specify the transmission rules across material boundaries.

Over the years, a variety of approaches have been developed to study interface problems in heat conduction through adjacent media. In certain cases, solution representations can be obtained using integral transform methods, such as Laplace and Fourier analysis, particularly under simple transmission conditions [1,2,4]. More recently, the celebrated *unified transform method* of Fokas has been applied successfully to the analysis of interface and contact problems, see e.g. [5-10,13]. This method provides a systematic framework for boundary and interface problems in linear PDEs, yielding explicit solution formulas under broad conditions. However, most existing studies have concentrated on homogeneous equations or restricted classes of interface conditions, while the general inhomogeneous interface problems in the heat conduction remain comparatively underdeveloped. This emphasizes the need for a broader setting that can accommodate interface problems for general heat conduction.

---

* corresponding author, e-mail address: chatziafrati@math.uoa.gr

In this paper, we investigate interface problems for the linear reaction-advection-diffusion (lRAD) equation in the presence of external sources, namely we study the following problem:

***Problem*** Solve

$$(1.1)_R \qquad \begin{cases} u_t^R = \sigma_R^2 u_{xx}^R + \gamma_R u_x^R + \delta_R u^R + f_R \\ u^R(x,0) = u_0^R(x) \end{cases}$$

for $u^R(x,t)$, with $x > 0$, $t > 0$, and

$$(1.1)_L \qquad \begin{cases} u_t^L = \sigma_L^2 u_{xx}^L + \gamma_L u_x^L + \delta_L u^L + f_L \\ u^L(x,0) = u_0^L(x) \end{cases}$$

for $u^L(x,t)$, with $x < 0$, $t > 0$, given the generic contact constraints:

$$(1.2) \qquad \begin{cases} \alpha_{11} u_x^R(0,t) + \alpha_{12} u^R(0,t) + \alpha_{13} u_x^L(0,t) + \alpha_{14} u^L(0,t) = I_1(t) \\ \alpha_{21} u_x^R(0,t) + \alpha_{22} u^R(0,t) + \alpha_{23} u_x^L(0,t) + \alpha_{24} u^L(0,t) = I_2(t) \end{cases}$$

for $t > 0$,

where the coefficients in (1.2) are real constants (further assumptions on these are made in the statement of our Proposition), the subscripts and superscripts $R$ and $L$ correspond to the right and left semi-infinite segments, respectively, the coefficients in equations (1.1) denote the diffusion, advection and reaction rates, respectively, in the two semi-infinite subdomains, while the nonhomogeneous terms represent external forcing. We assume that the initial conditions and the external forcing terms belong to the Schwartz class. More specifically, we have:

***Data*** We are given the initial values $u_0^R(x)$ ($x \ge 0$) and $u_0^L(x)$ ($x \le 0$), and the forcing functions $f_R(x,t)$ ($x \ge 0$, $t \ge 0$) and $f_L(x,t)$, ($x \le 0$, $t \ge 0$), satisfying the following assumptions:

$(1.3)_1$ $u_0^R \in \mathcal{S}([0,\infty))$, $u_0^L \in \mathcal{S}((-\infty,0])$,

$(1.3)_2$ $f_R(\cdot,t) \in \mathcal{S}([0,\infty))$, $f_L(\cdot,t) \in \mathcal{S}((-\infty,0])$, uniformly for $t$ in compact subsets of $[0,\infty)$.

***Parameters*** We make the following assumptions on the given parameters:

$(1.4)_1$ $\sigma_R > 0$, $\sigma_L > 0$,

$(1.4)_2$ either $\gamma_R \ge 0$ and $\delta_R \in \mathbb{R}$ or $\gamma_R < 0$ and $\delta_R < 0$,

$(1.4)_3$ either $\gamma_L \ge 0$ and $\delta_L \le 0$ or $\gamma_L < 0$ and $\delta_L \in \mathbb{R}$,

$(1.4)_4$ in any case we assume $\gamma_R^2 > 4\sigma_R^2 \delta_R$ and $\gamma_L^2 > 4\sigma_L^2 \delta_L$,

$(1.4)_5$ $\alpha_{jk} \in \mathbb{R}$.

The lRAD equation describes directed transport, diffusive spreading, and local reaction dynamics. The homogeneous lRAD equation is widely recognized as the advection (or convection)-diffusion equation, which models diffusive transport phenomena in physics, chemistry and biology [14-17]. In addition, the Black–Scholes equation, though originating in financial markets, shares the same mathematical structure as the lRAD equation: volatility induces diffusion, risk-free growth induces drift, and discounting acts as reaction [18]. Consequently, numerous analytical and numerical solutions have been developed for the one-, two-, and three-dimensional lRAD equations; see e.g. [19-21] and the references therein. In contrast, the analytical theory of inhomogeneous lRAD equations, particularly those involving general inhomogeneous interface problems, remains unexplored.

In this work, we derive explicit analytical representations of the solutions to Eqs. (1.1) incorporating both external forcing and interface transmission effects by means of the Fokas unified transform method (UTM). The UTM, developed as an extension of the inverse scattering transform, is a powerful tool for analyzing linear as well as nonlinear boundary value problems (see [22-34], in chronological order) with established analytical and numerical advantages compared to traditional transform methods. It has been successfully applied to not only evolution and elliptic PDEs but also, very recently, to problems of hyperbolic nature [35], while it has been extended to also address control problems [36,37]. The aim of this work is to extend these developments to the inhomogeneous lRAD equation with general interface conditions, providing explicit solution representations that account for both inhomogeneous forcing and the effects of interface transmission, thereby establishing a systematic framework for analyzing multi-domain lRAD systems.

Furthermore, this framework naturally connects with moving interface problems such as the classical Stefan problem, where the interface evolves dynamically according to energy balance laws; see [38] and literature cited therein. This connection highlights the broader applicability of our approach to phase transition phenomena. For the sake of completeness and the benefit of our readership, we should point out that further applications, recent advances and more bibliography related to the Fokas UTM methodology, which is a central tool in this study, can be found in e.g. [39-68], in alphabetical order.

The paper is organized as follows. In Section 2, we derive an integral representation for the solution of the lRAD equation for x > 0, while in Section 3 we present an analogous formula for the solution with x < 0, both using the UTM. As is well-known, a key component in the implementation of the UTM is the *global* relation, an integro-algebraic equation typically involving initial and boundary values, which in our case here can be used to characterize the unknown boundary values based on the prescribed interface conditions. In Section 4, we determine and eliminate the unknown boundary values by analyzing the global relations obtained for both x > 0 and x < 0, thereby deriving a closed-form analytical solution to the interface problem described in eqs. (1.1). We conclude with some remarks and the formulation of our Theorem.

## 2. Derivation of the solution for $x>0$

Setting

$$\hat{u}^R(\lambda,t)=(u^R\hat{)}(\lambda,t)=\int_{x=0}^{\infty}e^{-i\lambda x}u^R(x,t)dx \text{ and } \hat{f}_R(\lambda,t)=\int_{x=0}^{\infty}e^{-i\lambda x}f_R(x,t)dx\,, \text{ for } \operatorname{Im}\lambda\le 0\,,$$

we have

$$\frac{\partial \hat{u}^R(\lambda,t)}{\partial t}=\sigma_R^2\int_0^\infty e^{-i\lambda x}u_{xx}^R(x,t)dx+\gamma_R\int_0^\infty e^{-i\lambda x}u_x^R(x,t)dx+\delta_R\int_0^\infty e^{-i\lambda x}u^R(x,t)dx+\int_0^\infty e^{-i\lambda x}f_R(x,t)dx\,.$$

Integrating by parts we obtain

$$\frac{\partial \hat{u}^R(\lambda,t)}{\partial t}=-\sigma_R^2 g_1^R(t)+\sigma_R^2 i\lambda\int_0^\infty e^{-i\lambda x}u_x^R(x,t)dx-\gamma_R g_0^R(t)+\gamma_R i\lambda\int_0^\infty e^{-i\lambda x}u^R(x,t)dx+\delta_R\hat{u}^R(\lambda,t)+\hat{f}_R(\lambda,t)$$

and

$$\frac{\partial \hat{u}^R(\lambda,t)}{\partial t}=-\sigma_R^2 g_1^R(t)-\sigma_R^2 i\lambda g_0^R(t)-\sigma_R^2\lambda^2\hat{u}^R(\lambda,t)-\gamma_R g_0^R(t)+\gamma_R i\lambda\hat{u}^R(\lambda,t)+\delta_R\hat{u}^R(\lambda,t)+\hat{f}_R(\lambda,t)\,.$$

Thus,

(2.1) $$\frac{\partial \hat{u}^R(\lambda,t)}{\partial t}=-\omega_R(\lambda)\hat{u}^R(\lambda,t)-[\sigma_R^2 g_1^R(t)+\sigma_R^2 i\lambda g_0^R(t)+\gamma_R g_0^R(t)]+\hat{f}_R(\lambda,t) \text{ and } \hat{u}^R(\lambda,0)=\hat{u}_0^R(\lambda)\,,$$

where we have set

$$\omega_R(\lambda)=\sigma_R^2\lambda^2-i\gamma_R\lambda-\delta_R\,.$$

Solving the ODE (2.1), we find the global relation

(2.2) $$\hat{u}^R(\lambda,t)=e^{-\omega_R(\lambda)t}\hat{u}_0^R(\lambda)-e^{-\omega_R(\lambda)t}[\sigma_R^2\tilde{g}_1^R(\omega_R,t)+(\sigma_R^2 i\lambda+\gamma_R)\tilde{g}_0^R(\omega_R,t)]+e^{-\omega_R(\lambda)t}\tilde{\hat{f}}_R(\lambda,\omega_R,t)\,,\quad \operatorname{Im}\lambda\le 0\,,$$

where $g_1^R(t):=u_x^R(0,t)$, $g_0^R(t):=u_x^R(0,t)$, and

$$\tilde{g}_1^R(\omega_R,t)=\int_{\tau=0}^{t}e^{\omega_R(\lambda)\tau}g_1^R(\tau)d\tau\,,\ \tilde{g}_0^R(\omega_R,t)=\int_{\tau=0}^{t}e^{\omega_R(\lambda)\tau}g_0^R(\tau)d\tau\ \text{ and }\ \tilde{\hat{f}}_R(\lambda,\omega_R,t)=\int_{\tau=0}^{t}e^{\omega_R(\lambda)\tau}\hat{f}_R(\lambda,\tau)d\tau\,.$$

*Comment* An alternative way to obtain the general relation (2.2) is based on Green's theorem and the equation written in divergence form, i.e.,

$$\left(e^{i\lambda x-\omega_R(\lambda)t}u\right)_t=\left\{e^{i\lambda x-\omega_R(\lambda)t}[(i\sigma_R^2\lambda+\gamma_R)u^R+\sigma_R^2u_x^R]\right\}_x+e^{i\lambda x-\omega_R(\lambda)t}f_R(x,t)\,.$$

The above form of the equation is directly connected to the Lax pair associated with the homogeneous problem, i.e., when $f_R\equiv 0$:

$$\mathcal{L}_x-i\lambda\mathcal{L}=u^R\,,\ \mathcal{L}_t-\omega_R(\lambda)\mathcal{L}=(i\sigma_R^2\lambda+\gamma_R)u^R+\sigma_R^2u_x^R\,.$$

By Fourier inversion formula,

(2.3) $$u^R(x,t)=\frac{1}{2\pi}\int_{-\infty}^{\infty}e^{i\lambda x-\omega_R(\lambda)t}\hat{u}_0^R(\lambda)d\lambda-\frac{1}{2\pi}\int_{-\infty}^{\infty}e^{i\lambda x-\omega_R(\lambda)t}[\sigma_R^2\tilde{g}_1^R(\omega_R,t)+(\sigma_R^2i\lambda+\gamma_R)\tilde{g}_0^R(\omega_R,t)]d\lambda$$
$$+\frac{1}{2\pi}\int_{-\infty}^{\infty}e^{i\lambda x-\omega_R(\lambda)t}\tilde{\hat{f}}_R(\lambda,\omega_R,t)d\lambda\,,\ \text{for } x>0\,,\ t>0\,.$$

Deforming the contour, we write equation (2.3) in the following way:

(2.4) $$u^R(x,t)=\frac{1}{2\pi}\int_{-\infty}^{\infty}e^{i\lambda x-\omega_R(\lambda)t}\hat{u}_0^R(\lambda)d\lambda-\frac{1}{2\pi}\int_{\lambda\in C_R^+}e^{i\lambda x-\omega_R(\lambda)t}[\sigma_R^2\tilde{g}_1^R(\omega_R,t)+(\sigma_R^2i\lambda+\gamma_R)\tilde{g}_0^R(\omega_R,t)]d\lambda$$
$$+\frac{1}{2\pi}\int_{-\infty}^{\infty}e^{i\lambda x-\omega_R(\lambda)t}\tilde{\hat{f}}_R(\lambda,\omega_R,t)d\lambda\,,\ \text{for } x>0\,,\ t>0\,,$$

where $C_R^+$ is the contour (fig 1)

$$C_R^+=\{\lambda\in\mathbb{C}:\operatorname{Im}\lambda\geq 0\ \textit{and}\ \operatorname{Re}[\omega_R(\lambda)]=\operatorname{Re}[\sigma_R^2\lambda^2-i\gamma_R\lambda-\delta_R]=0\}\,.$$

Writing $\lambda=\xi+i\eta$, with $\xi,\eta\in\mathbb{R}$, the contour $C_R^+$ is the graph of the function

$$\eta=\frac{\gamma_R}{2\sigma_R^2}+\sqrt{\xi^2+\frac{\gamma_R^2-4\sigma_R^2\delta_R}{4\sigma_R^4}}\,,\ \xi\in\mathbb{R}\,,$$

in the $\xi\eta-$plane.

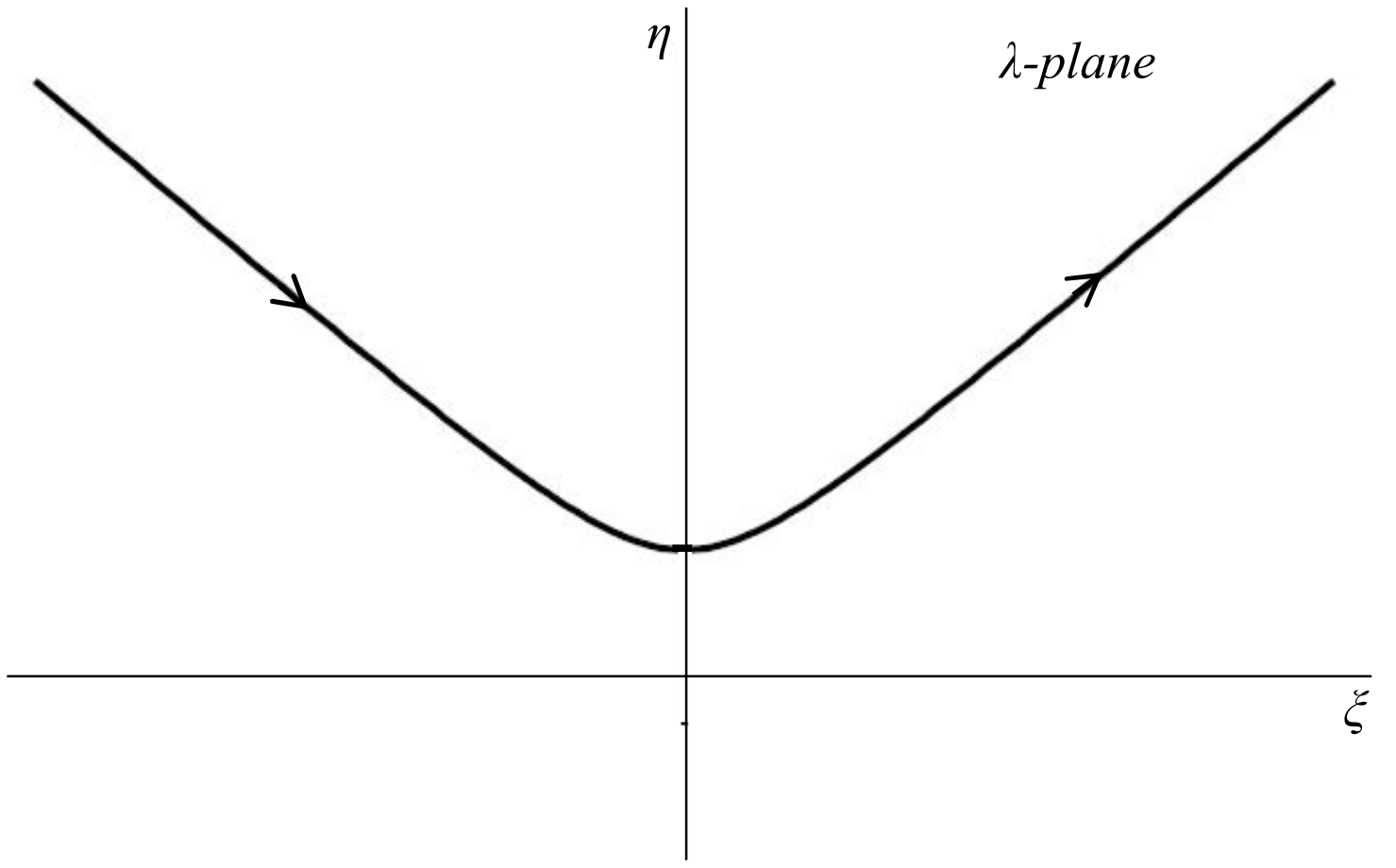

**Fig 1** *The contour* $C_R^+$. *By the assumptions (1.4),* $\eta>0$, $\forall(\xi,\eta)\in C_R^+$.

Also, (2.4) can be written as follows: For every $T>0$,

(2.5) $u^R(x,t)=\frac{1}{2\pi}\int_{-\infty}^{\infty}e^{i\lambda x-\omega_R(\lambda)t}\hat{u}_0^R(\lambda)d\lambda-\frac{1}{2\pi}\int_{\lambda\in C_R^+}e^{i\lambda x-\omega_R(\lambda)t}[\sigma_R^2\tilde{g}_1^R(\omega_R,T)+(\sigma_R^2 i\lambda+\gamma_R)\tilde{g}_0^R(\omega_R,T)]d\lambda$

$$+\frac{1}{2\pi}\int_{-\infty}^{\infty}e^{i\lambda x-\omega_R(\lambda)t}\tilde{\hat{f}}_R(\lambda,\omega_R,t)d\lambda\text{, for } x>0\text{, } 0<t<T.$$

Indeed, this is a consequence of Cauchy's theorem and Jordan's lemma.

Since $\omega_R(-\lambda+\frac{i\gamma_R}{\sigma_R^2})=\omega_R(\lambda)$, (2.2) gives

(2.6) $e^{\omega_R(\lambda)t}\hat{u}^R(-\lambda+\frac{i\gamma_R}{\sigma_R^2},t)=\hat{u}_0^R(-\lambda+\frac{i\gamma_R}{\sigma_R^2})-\sigma_R^2[\tilde{g}_1^R(\omega_R,t)-i\lambda\tilde{g}_0^R(\omega_R,t)]+\tilde{\hat{f}}_R(-\lambda+\frac{i\gamma_R}{\sigma_R^2},\omega_R,t)$, $\operatorname{Im}\lambda\ge\frac{\gamma_R}{\sigma_R^2}$, $t>0$.

Let us write the above equation also in the form

(2.7) $e^{\omega_R(\lambda)T}\hat{u}^R(-\lambda+\frac{i\gamma_R}{\sigma_R^2},T)=\hat{u}_0^R(-\lambda+\frac{i\gamma_R}{\sigma_R^2})-\sigma_R^2[\tilde{g}_1^R(\omega_R,T)-i\lambda\tilde{g}_0^R(\omega_R,T)]$

$$+\tilde{\hat{f}}_R(-\lambda+\frac{i\gamma_R}{\sigma_R^2},\omega_R,T),\ \operatorname{Im}\lambda\ge\frac{\gamma_R}{\sigma_R^2},T>0.$$

Next we solve the equation $\omega_R(\lambda)=\sigma_R^2\lambda^2-i\gamma_R\lambda-\delta_R=\mu^2$ for $\lambda$, by considering the function

$$\lambda=\varphi_R(\mu):=\frac{i\gamma_R}{2\sigma_R^2}+\frac{1}{\sigma_R}\sqrt[R]{\mu^2-\varepsilon_R^2}\text{ , } \mu\in\mathbb{C}-[-\varepsilon_R,\varepsilon_R],$$

where $\varepsilon_R=\sqrt{(\gamma_R^2-4\sigma_R^2\delta_R)/(4\sigma_R^2)}>0$ and the $R-$*square root* function $\sqrt[R]{\mu^2-\varepsilon_R^2}$ is defined as follows:
If

$$\mu-\varepsilon_R=|\mu-\varepsilon_R|e^{i\theta_R^+(\mu)}\text{ and }\mu-(-\varepsilon_R)=|\mu-(-\varepsilon_R)|e^{i\theta_R^-(\mu)}\text{, }0\le\theta_R^+(\mu),\theta_R^-(\mu)<2\pi\text{ (fig 2),}$$

then

$$\sqrt[R]{\mu^2-\varepsilon_R^2}:=\sqrt{|\mu^2-\varepsilon_R^2|}\exp\left(i\frac{\theta_R^+(\mu)+\theta_R^-(\mu)}{2}\right).$$

*Note.* $R$, in the symbol $\sqrt[R]{\mu^2-\varepsilon_R^2}$, stands for "*RIGHT*". Also, with this choice of the angles $\theta_R^+(\mu)$ and $\theta_R^-(\mu)$, $\sqrt[R]{(-\mu)^2-\varepsilon_R^2}=-\sqrt[R]{\mu^2-\varepsilon_R^2}$, despite the presence of the term "$\mu^2$".

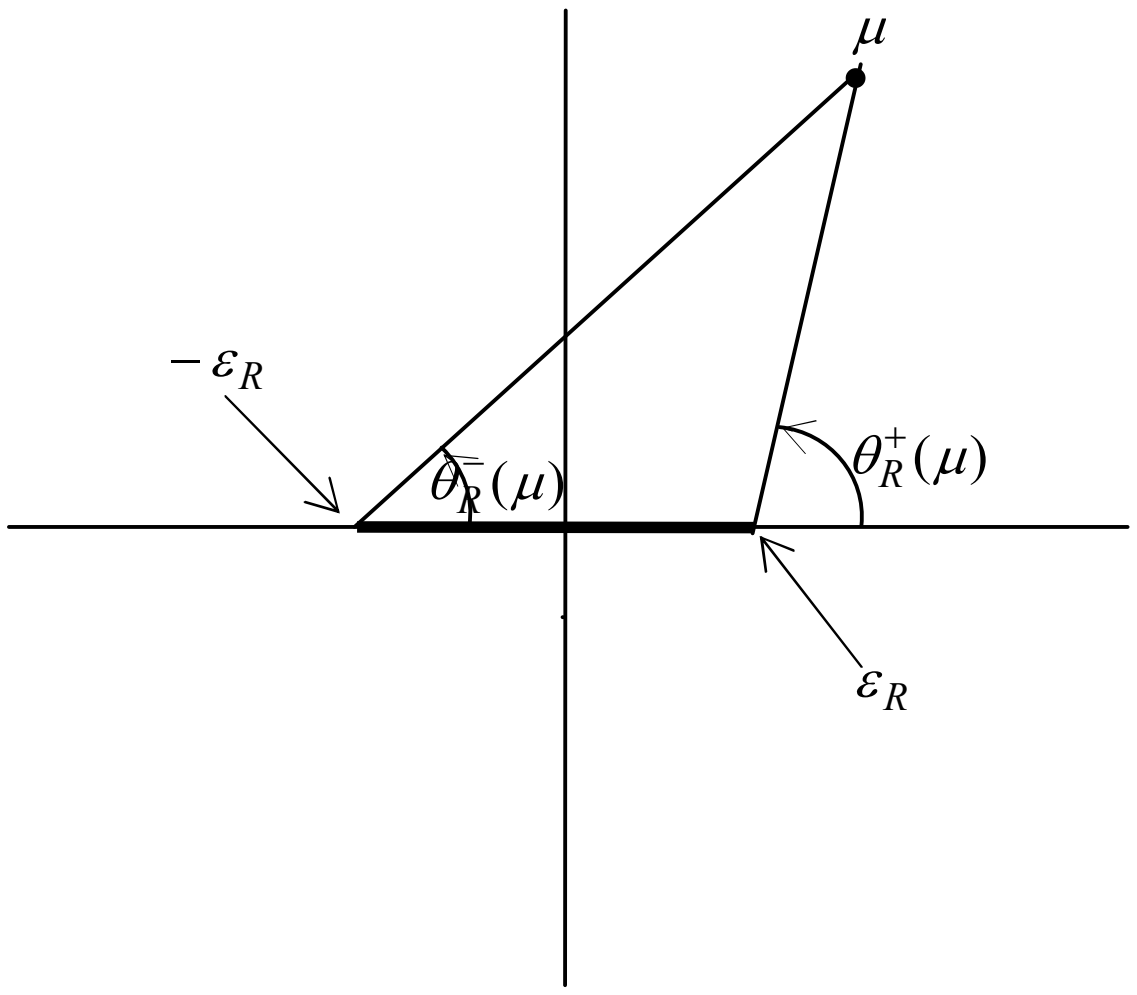

**Fig.2** *The angles* $\theta_R^+(\mu)$ *and* $\theta_R^-(\mu)$ *in the definition of the function* $\sqrt[R]{\mu^2-\varepsilon_R^2}$.

We point out that the function $\varphi_R(\mu)$ is analytic (holomorphic) for $\mu \in \mathbb{C} - [-\varepsilon_R, \varepsilon_R]$ and has continuous extension to the half-plane $\{\mu \in \mathbb{C} : \operatorname{Im}\mu \geq 0\}$. Moreover this extension is $C^\infty$ in $\{\mu \in \mathbb{C} : \operatorname{Im}\mu \geq 0\} - \{-\varepsilon_R, \varepsilon_R\}$. Let $\Delta^+$ be the set $\Delta^+ = \{\mu \in \mathbb{C} : \operatorname{Im}\mu \geq 0 \ \textit{and} \ \operatorname{Re}[\mu^2] \leq 0\}$ and $\Gamma^+$ be its boundary, i.e., the contour

$$\Gamma^+ = \partial\Delta^+ = \{\mu \in \mathbb{C} : \operatorname{Im}\mu \geq 0 \ \textit{and} \ \operatorname{Re}[\mu^2] = 0\}. \text{ (Fig.3)}$$

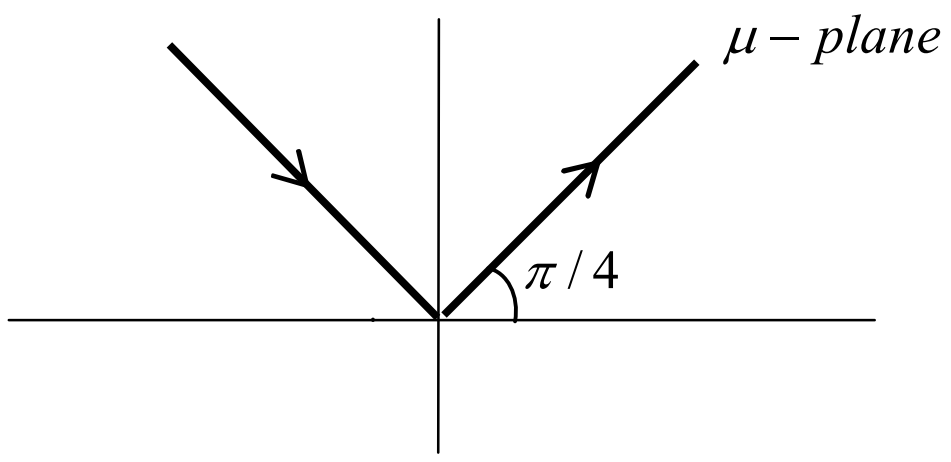

**Fig 3** *The contour* $\Gamma^+$

Changing the variable in the integrals over $C_R^+$, in (2.5), from $\lambda \in C_R^+$ to $\mu \in \Gamma^+$ by setting $\lambda = \varphi_R(\mu)$, we obtain

$$(2.8)\, u^R(x,t) = \frac{1}{2\pi}\int_{-\infty}^{\infty} e^{i\lambda x - \omega_R(\lambda)t}\hat{u}_0^R(\lambda)d\lambda$$
$$-\frac{1}{2\pi}\int_{\mu\in\Gamma^+} e^{i\varphi_R(\mu)x - \mu^2 t}\{\sigma_R^2\tilde{g}_1^R(\mu^2,T) + [\sigma_R^2 i\varphi_R(\mu) + \gamma_R]\tilde{g}_0^R(\mu^2,T)\}\frac{d\varphi_R(\mu)}{d\mu}d\mu$$
$$+\frac{1}{2\pi}\int_{-\infty}^{\infty} e^{i\lambda x - \omega_R(\lambda)t}\tilde{\hat{f}}_R(\lambda,\omega_R,t)d\lambda, \text{ for } x > 0,\ 0 < t < T.$$

*Remark* Let as keep in mind that

$$\frac{d\lambda}{d\mu} = \frac{d\varphi_R(\mu)}{d\mu} = \frac{1}{\sigma_R}\frac{\mu}{{}^R\!\sqrt{\mu^2 - \varepsilon_R^2}},\quad \lim_{\substack{\mu\to 0\\ \operatorname{Im}\mu>0}}\frac{d\varphi_R(\mu)}{d\mu} = 0 \text{ and } \frac{d\varphi_R(\mu)}{d\mu} = O(1), \text{ as } \mu \to \infty.$$

Also, (2.6) gives

$$(2.9)\quad e^{\mu^2 T}\hat{u}^R(-\varphi_R(\mu) + \tfrac{i\gamma_R}{\sigma_R^2}, T) = \hat{u}_0^R(-\varphi_R(\mu) + \tfrac{i\gamma_R}{\sigma_R^2}) - \sigma_R^2[\tilde{g}_1^R(\mu^2,T) - i\varphi_R(\mu)\tilde{g}_0^R(\mu^2,T)]$$
$$+\tilde{\hat{f}}_R(-\varphi_R(\mu) + \tfrac{i\gamma_R}{\sigma_R^2}, \mu^2, T),\ \operatorname{Im}\varphi_R(\mu) \geq \tfrac{\gamma_R}{\sigma_R^2}.$$

In particular (2.9) holds for $\mu \in \Delta^+$.

## 3. Derivation of the solution for $x < 0$

As in *Section 2*, setting

$$\hat{u}^L(\lambda,t) = (u^L\hat{)}(\lambda,t) = \int_{x=-\infty}^{0} e^{-i\lambda x}u^L(x,t)dx \text{ and } \hat{f}_L(\lambda,t) = \int_{x=0}^{\infty} e^{-i\lambda x}f_L(x,t)dx, \text{ for } \operatorname{Im}\lambda \geq 0,$$

we have

$$\frac{\partial \hat{u}^L(\lambda,t)}{\partial t} = \sigma_L^2\int_{-\infty}^{0} e^{-i\lambda x}u_{xx}^L(x,t)dx + \gamma_L\int_{-\infty}^{0} e^{-i\lambda x}u_x^L(x,t)dx + \delta_L\int_{-\infty}^{0} e^{-i\lambda x}u^L(x,t)dx + \int_{-\infty}^{0} e^{-i\lambda x}f_L(x,t)dx.$$

Integrating by parts as before, we obtain

$$(3.1)\qquad \frac{\partial \hat{u}^L(\lambda,t)}{\partial t} = -\omega_L(\lambda)\hat{u}^L(\lambda,t) + [\sigma_L^2 g_1^L(t) + \sigma_L^2 i\lambda g_0^L(t) + \gamma_L g_0^L(t)] + \hat{f}_L(\lambda,t) \text{ and } \hat{u}^L(\lambda,0) = \hat{u}_0^L(\lambda),$$

where we have set

$$\omega_L(\lambda)=\sigma_L^2\lambda^2-i\gamma_L\lambda-\delta_L .$$

Solving the ODE (3.1), we find

$$(3.2)\quad \hat{u}^L(\lambda,t)=e^{-\omega_L(\lambda)t}\hat{u}_0^L(\lambda)+e^{-\omega_R(\lambda)t}[\sigma_L^2\tilde{g}_1^L(\omega_L,t)+(\sigma_L^2 i\lambda+\gamma_L)\tilde{g}_0^L(\omega_L,t)]+e^{-\omega_L(\lambda)t}\tilde{\hat{f}}_L(\lambda,\omega_L,t),\quad \operatorname{Im}\lambda\ge 0,$$

Where $g_1^L(t):=u_x^L(0,t)$, $g_0^L(t):=u_x^L(0,t)$, and

$$\tilde{g}_1^L(\omega_L,t)=\int_{\tau=0}^{t}e^{-\omega_L(\lambda)\tau}g_1^L(\tau)d\tau\ ,\ \tilde{g}_0^L(\omega_L,t)=\int_{\tau=0}^{t}e^{-\omega_L(\lambda)\tau}g_0^L(\tau)d\tau\ \text{ and }\ \tilde{\hat{f}}_L(\lambda,\omega_L,t)=\int_{\tau=0}^{t}e^{-\omega_L(\lambda)\tau}\hat{f}_L(\lambda,\tau)d\tau .$$

By Fourier inversion formula, (3.2) gives

$$(3.3)\ u^L(x,t)=\frac{1}{2\pi}\int_{-\infty}^{\infty}e^{i\lambda x-\omega_L(\lambda)t}\hat{u}_0^L(\lambda)d\lambda+\frac{1}{2\pi}\int_{-\infty}^{\infty}e^{i\lambda x-\omega_L(\lambda)t}[\sigma_L^2\tilde{g}_1^L(\omega_L,t)+(\sigma_L^2 i\lambda+\gamma_L)\tilde{g}_0^L(\omega_L,t)]d\lambda$$

$$+\frac{1}{2\pi}\int_{-\infty}^{\infty}e^{i\lambda x-\omega_L(\lambda)t}\tilde{\hat{f}}_L(\lambda,\omega_L,t)d\lambda,\ \text{for } x<0,\ t>0,$$

Also (3.3) can be written as follows:

$$(3.4)\ u^L(x,t)=\frac{1}{2\pi}\int_{-\infty}^{\infty}e^{i\lambda x-\omega_L(\lambda)t}\hat{u}_0^L(\lambda)d\lambda+\frac{1}{2\pi}\int_{\lambda\in C_L^-}e^{i\lambda x-\omega_L(\lambda)t}[\sigma_L^2\tilde{g}_1^L(\omega_L,t)+(\sigma_L^2 i\lambda+\gamma_L)\tilde{g}_0^L(\omega_L,t)]d\lambda$$

$$+\frac{1}{2\pi}\int_{-\infty}^{\infty}e^{i\lambda x-\omega_L(\lambda)t}\tilde{\hat{f}}_L(\lambda,\omega_L,t)d\lambda,\ \text{for } x<0,\ t>0,$$

where $C_L^-$ is the contour (fig 4)

$$C_L^-=\{\lambda\in\mathbb{C}:\operatorname{Im}\lambda\le 0\ \mathit{and}\ \operatorname{Re}[\omega_L(\lambda)]=\operatorname{Re}[\sigma_L^2\lambda^2-i\gamma_L\lambda-\delta_L]=0\}.$$

Writing $\lambda=\xi+i\eta$, with $\xi,\eta\in\mathbb{R}$, the contour $C_L^-$ is the graph of the function

$$\eta=\frac{\gamma_L}{2\sigma_L^2}-\sqrt{\xi^2+\frac{\gamma_L^2-4\sigma_L^2\delta_L}{4\sigma_L^4}},\ \xi\in\mathbb{R},$$

in the $\xi\eta-$plane.

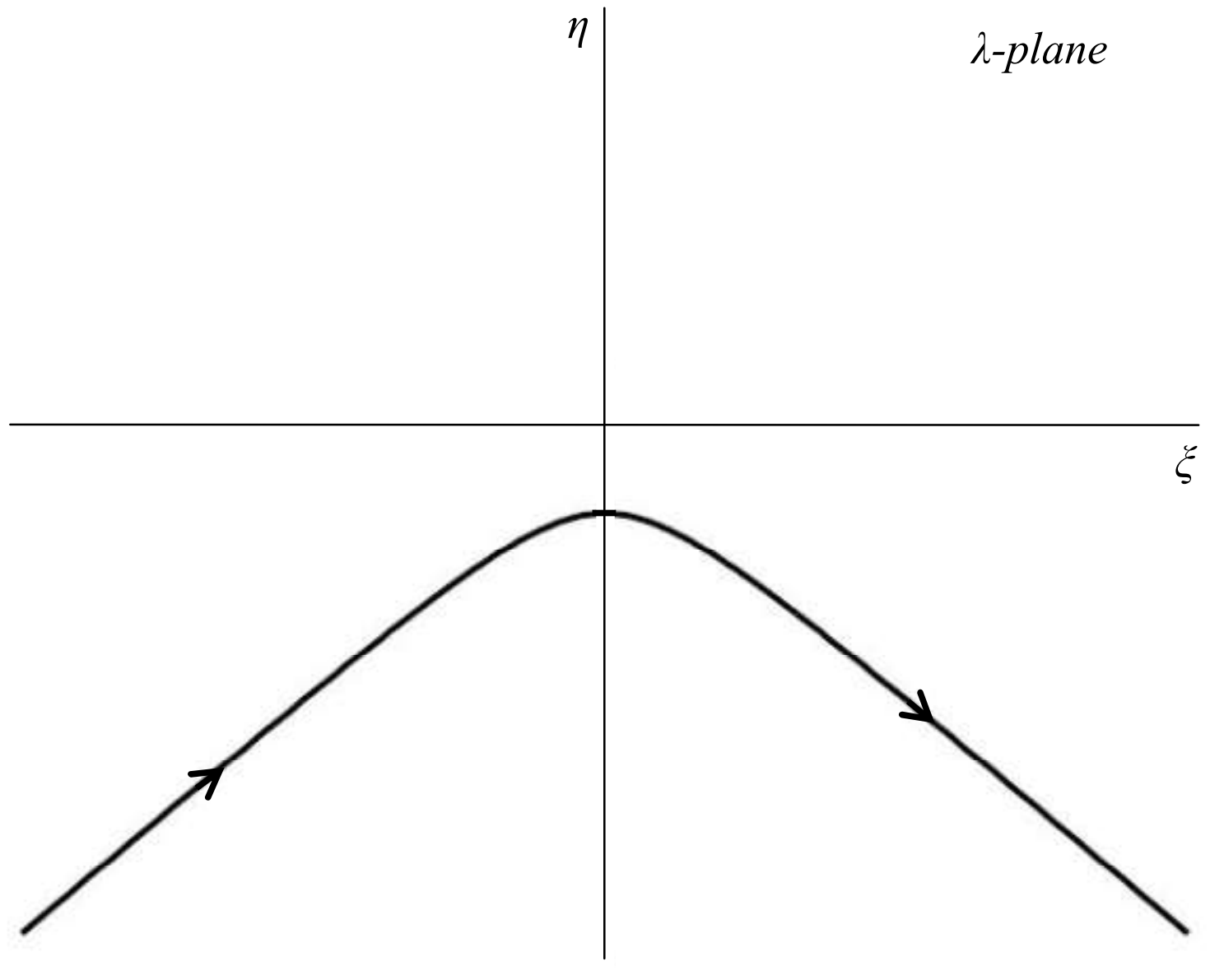

**Fig 4** *The contour* $C_L^-$. *By the assumptions (1.4),* $\eta<0$, $\forall(\xi,\eta)\in C_L^-$.

Also, for every $T>0$, (3.4) can be written as follows:

$$(3.5)\quad u^L(x,t)=\frac{1}{2\pi}\int_{-\infty}^{\infty}e^{i\lambda x-\omega_L(\lambda)t}\hat{u}_0^L(\lambda)d\lambda+\frac{1}{2\pi}\int_{\lambda\in C_L^-}e^{i\lambda x-\omega_L(\lambda)t}[\sigma_L^2\tilde{g}_1^L(\omega_L,T)+(\sigma_L^2 i\lambda-\gamma_L)\tilde{g}_0^L(\omega_L,T)]d\lambda$$

$$+\frac{1}{2\pi}\int_{-\infty}^{\infty}e^{i\lambda x-\omega_L(\lambda)t}\tilde{\hat{f}}_L(\lambda,\omega_L,t)d\lambda\text{, for } x<0\text{, } 0<t<T.$$

Since $\omega_L(-\lambda+\frac{i\gamma_L}{\sigma_L^2})=\omega_L(\lambda)$, (3.2) gives

$$(3.6)\quad e^{\omega_L(\lambda)t}\hat{u}^L(-\lambda+\tfrac{i\gamma_L}{\sigma_L^2},t)=\hat{u}_0^L(-\lambda+\tfrac{i\gamma_L}{\sigma_L^2})+\sigma_L^2[\tilde{g}_1^L(\omega_L,t)-i\lambda\tilde{g}_0^L(\omega_L,t)]+\tilde{\hat{f}}_L(-\lambda+\tfrac{i\gamma_L}{\sigma_L^2},\omega_L,t),\ \operatorname{Im}\lambda\le\tfrac{\gamma_L}{\sigma_L^2},\ t>0.$$

Also, we have

$$(3.7)\quad e^{\omega_L(\lambda)T}\hat{u}^L(-\lambda+\tfrac{i\gamma_L}{\sigma_L^2},T)=\hat{u}_0^L(-\lambda+\tfrac{i\gamma_L}{\sigma_L^2})+\sigma_L^2[\tilde{g}_1^L(\omega_L,T)-i\lambda\tilde{g}_0^L(\omega_L,T)]$$

$$+\tilde{\hat{f}}_L(-\lambda+\tfrac{i\gamma_L}{\sigma_L^2},\omega_L,T),\ \operatorname{Im}\lambda\le\tfrac{\gamma_L}{\sigma_L^2},\ T>0.$$

Next we solve the equation $\omega_L(\lambda)=\sigma_L^2\lambda^2-i\gamma_L\lambda-\delta_L=\mu^2$ for $\lambda$, by considering the function

$$\lambda=\varphi_L(\mu):=\frac{i\gamma_L}{2\sigma_L^2}+\frac{1}{\sigma_L}\sqrt[L]{\mu^2-\varepsilon_L^2},\ \mu\in\mathbb{C}-[-\varepsilon_L,\varepsilon_L],$$

where $\varepsilon_L=\sqrt{(\gamma_L^2-4\sigma_L^2\delta_L)/(4\sigma_L^2)}>0$ and the $L-$*square root* function $\sqrt[L]{\mu^2-\varepsilon_L^2}$ is defined, *this time*, as follows: If

$$\mu-\varepsilon_L=|\mu-\varepsilon_L|e^{i\theta_L^+(\mu)}\text{ and }\mu-(-\varepsilon_L)=|\mu-(-\varepsilon_L)|e^{i\theta_L^-(\mu)},\ -2\pi<\theta_L^+(\mu),\theta_L^-(\mu)\le 0\text{ (fig 5),}$$

then

$$\sqrt[L]{\mu^2-\varepsilon_L^2}:=\sqrt{|\mu^2-\varepsilon_L^2|}\exp\left(i\frac{\theta_L^+(\mu)+\theta_L^-(\mu)}{2}\right).$$

*Note*. $L$, in the symbol $\sqrt[L]{\mu^2-\varepsilon_L^2}$, stands for "*LEFT*". Also, with this choice of the angles $\theta_L^+(\mu)$ and $\theta_L^-(\mu)$,

$$\sqrt[L]{(-\mu)^2-\varepsilon_L^2}=-\sqrt[L]{\mu^2-\varepsilon_L^2},$$

despite the presence of the term "$\mu^2$". Thus,

$$\varphi_L(-\mu)=\frac{i\gamma_L}{2\sigma_L^2}-\frac{1}{\sigma_L}\sqrt[L]{\mu^2-\varepsilon_L^2}.$$

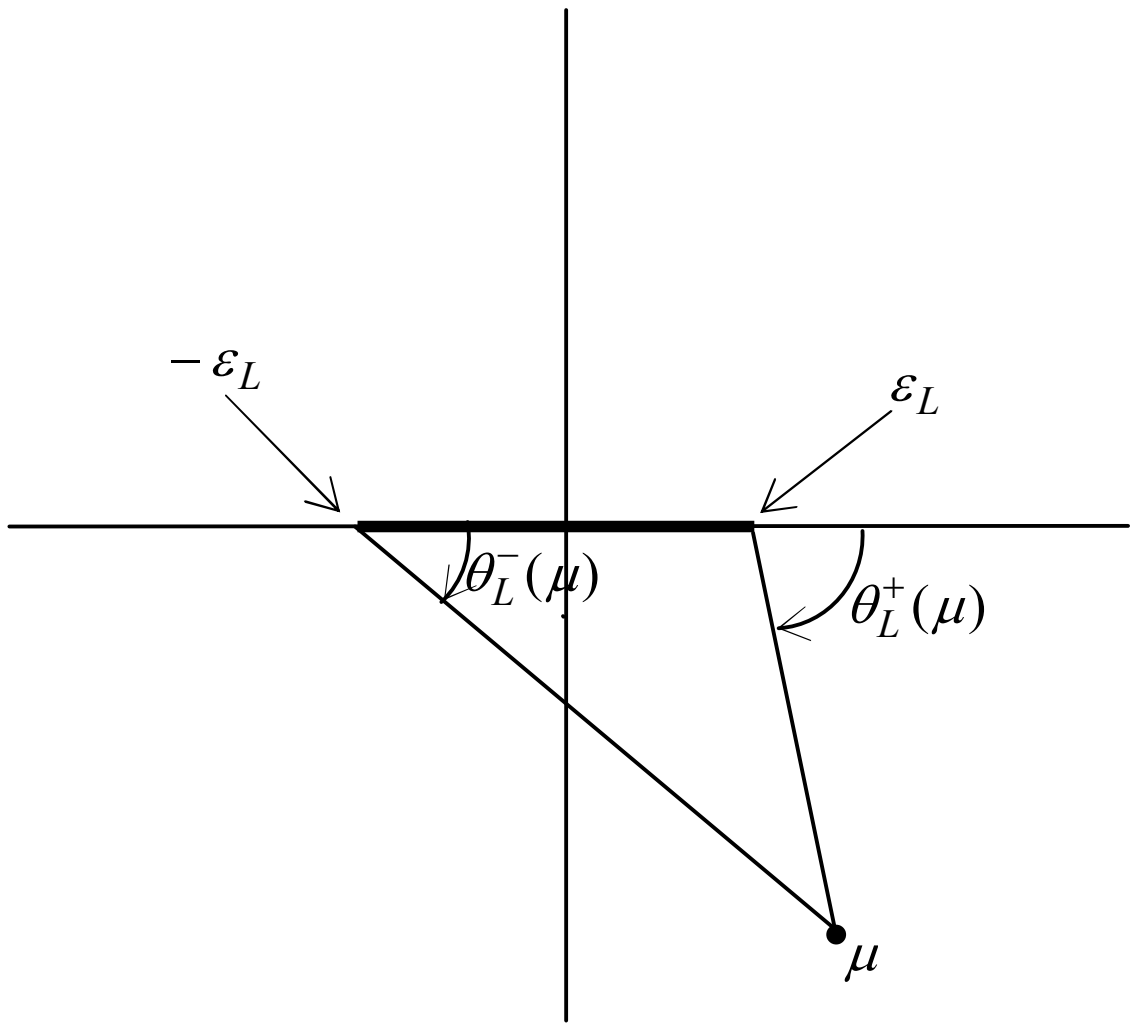

**Fig.5** *The angles* $\theta_L^+(\mu)$ *and* $\theta_L^-(\mu)$ *in the definition of the function* $\sqrt[L]{\mu^2-\varepsilon_L^2}$ .

The function $\varphi_L(\mu)$ is analytic for $\mu\in\mathbb{C}-[-\varepsilon_L,\varepsilon_L]$ and has continuous extension to the half-plane $\{\mu\in\mathbb{C}:\operatorname{Im}\mu\le 0\}$. Moreover this extension is $C^\infty$ in $\{\mu\in\mathbb{C}:\operatorname{Im}\mu\le 0\}-\{-\varepsilon_L,\varepsilon_L\}$.

Let $\Gamma^-$ be the contour (fig 6) $\Gamma^-=\{\mu\in\mathbb{C}:\operatorname{Im}\mu\le 0 \ \textit{and}\ \operatorname{Re}[\mu^2]=0\}$.

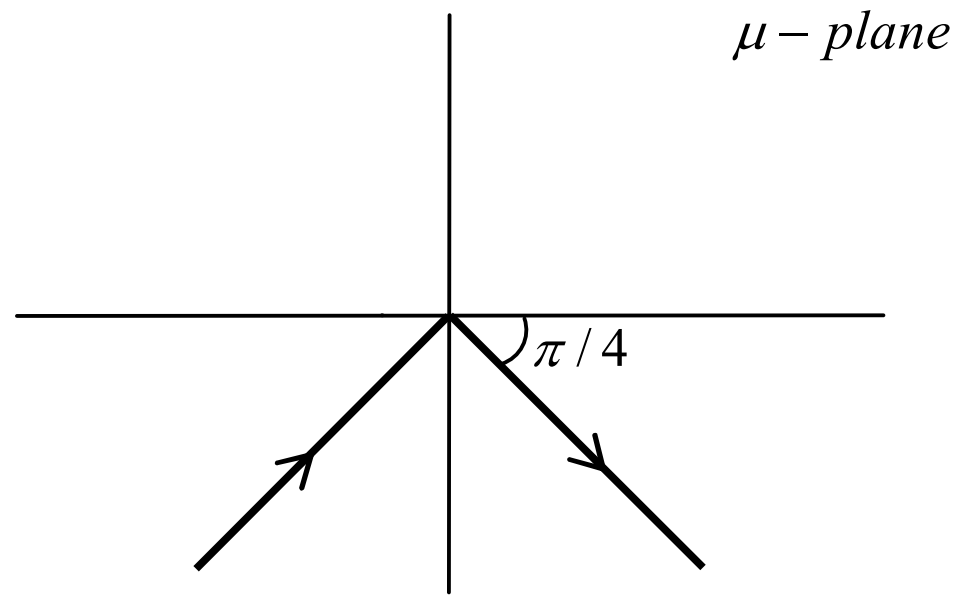

**Fig.6** *The contour* $\Gamma^-$

Changing the variable in the integrals over $C_L^-$, in (3.5), from $\lambda\in C_L^-$ to $\mu\in\Gamma^-$ by setting $\lambda=\varphi_L(\mu)$, we obtain

$$u^L(x,t)=\frac{1}{2\pi}\int_{-\infty}^{\infty}e^{i\lambda x-\omega_L(\lambda)t}\hat{u}_0^L(\lambda)d\lambda+\frac{\sigma_L^2}{2\pi}\int_{\mu\in\Gamma^-}e^{i\varphi_L(\mu)x-\mu^2t}\{\tilde{g}_1^L(\mu^2,T)-[i\varphi_L(\mu)+\gamma_L]\tilde{g}_0^L(\mu^2,T)\}\frac{d\varphi_L(\mu)}{d\mu}d\mu$$

$$+\frac{1}{2\pi}\int_{-\infty}^{\infty}e^{i\lambda x-\omega_L(\lambda)t}\tilde{\hat{f}}_L(\lambda,\omega_L,t)d\lambda,\ \text{for } x<0,\ 0<t<T.$$

We can write the above equation as follows:

$$(3.8)\quad u^L(x,t)=\frac{1}{2\pi}\int_{-\infty}^{\infty}e^{i\lambda x-\omega_L(\lambda)t}\hat{u}_0^L(\lambda)d\lambda$$

$$+\frac{\sigma_L^2}{2\pi}\int_{\mu\in\Gamma^+}e^{i\varphi_L(-\mu)x-\mu^2t}\{\tilde{g}_1^L(\mu^2,T)+[i\varphi_L(-\mu)+\gamma_L]\tilde{g}_0^L(\mu^2,T)\}\left\{-\frac{d[\varphi_L(-\mu)]}{d\mu}\right\}d\mu$$

$$+\frac{1}{2\pi}\int_{-\infty}^{\infty}e^{i\lambda x-\omega_L(\lambda)t}\tilde{\hat{f}}_L(\lambda,\omega_L,t)d\lambda,\ \text{for } x<0,\ 0<t<T.$$

Moreover, (3.7) gives

$$(3.9)\quad e^{\mu^2T}\hat{u}^L(-\varphi_L(\mu)+\tfrac{i\gamma_L}{\sigma_L^2},T)=\hat{u}_0^L(-\varphi_L(\mu)+\tfrac{i\gamma_L}{\sigma_L^2})-\sigma_L^2[\tilde{g}_1^L(\mu^2,T)-i\varphi_L(\mu)\tilde{g}_0^L(\mu^2,T)]$$

$$+\tilde{\hat{f}}_L(-\varphi_L(\mu)+\frac{i\gamma_L}{\sigma_L^2},\mu^2,T),\ \operatorname{Im}\varphi_L(\mu)\le\frac{\gamma_L}{\sigma_L^2}.$$

In particular the above equation holds for $\mu\in\Gamma^-$. Now, if $\mu\in\Gamma^+$ then $-\mu\in\Gamma^-$, and therefore (3.9) gives

$$(3.10)\ e^{\mu^2T}\hat{u}^L(-\varphi_L(-\mu)+\frac{i\gamma_L}{\sigma_L^2},T)=\hat{u}_0^L(-\varphi_L(-\mu)+\frac{i\gamma_L}{\sigma_L^2})-\sigma_L^2[\tilde{g}_1^L(\mu^2,T)-i\varphi_L(-\mu)\tilde{g}_0^L(\mu^2,T)]$$
$$+\tilde{\hat{f}}_L(-\varphi_L(-\mu)+\frac{i\gamma_L}{\sigma_L^2},\mu^2,T)\ (\text{for }\mu\in\Delta^+).$$

## 4. Computing the integrals containing unknown quantities

Equations (2.9) and (3.10) can be written in the following way:

$$(4.1)_R\quad \tilde{g}_1^R(\mu^2,T)-i\varphi_R(\mu)\tilde{g}_0^R(\mu^2,T)=J_R^*(\mu,T)+J_R(\mu,T),\ \mu\in\Delta^+,$$

$$(4.1)_L\quad \tilde{g}_1^L(\mu^2,T)-i\varphi_L(-\mu)\tilde{g}_0^L(\mu^2,T)=J_L^*(\mu,T)+J_L(\mu,T),\ \mu\in\Delta^+,$$

where we have set

$$J_R^*(\mu,T)=-\frac{1}{\sigma_R^2}e^{\mu^2T}\hat{u}^R(-\varphi_R(\mu)+\frac{i\gamma_R}{\sigma_R^2},T)$$

$$J_L^*(\mu,T)=-\frac{1}{\sigma_L^2}e^{\mu^2T}\hat{u}^L(-\varphi_L(-\mu)+\frac{i\gamma_L}{\sigma_L^2},T)$$

$$J_R(\mu,T)=\frac{1}{\sigma_R^2}\hat{u}_0^R(-\varphi_R(\mu)+\frac{i\gamma_R}{\sigma_R^2})+\frac{1}{\sigma_R^2}\tilde{\hat{f}}_R(-\varphi_R(\mu)+\frac{i\gamma_R}{\sigma_R^2},\mu^2,T)$$

$$J_L(\mu,T)=\frac{1}{\sigma_L^2}\hat{u}_0^L(-\varphi_L(-\mu)+\frac{i\gamma_L}{\sigma_L^2})+\frac{1}{\sigma_L^2}\tilde{\hat{f}}_L(-\varphi_L(-\mu)+\frac{i\gamma_L}{\sigma_L^2},\mu^2,T).$$

On the other hand, the interface conditions (1.2) give the following equations:

$$(4.2)_1\quad \alpha_{11}\tilde{g}_1^R(\mu^2,T)+\alpha_{12}\tilde{g}_0^R(\mu^2,T)+\alpha_{13}\tilde{g}_1^L(\mu^2,T)+\alpha_{14}\tilde{g}_0^L(\mu^2,T)=\tilde{I}_1(\mu^2,T),\ \mu\in\mathbb{C},$$

$$(4.2)_2\quad \alpha_{21}\tilde{g}_1^R(\mu^2,T)+\alpha_{22}\tilde{g}_0^R(\mu^2,T)+\alpha_{23}\tilde{g}_1^L(\mu^2,T)+\alpha_{24}\tilde{g}_0^L(\mu^2,T)=\tilde{I}_2(\mu^2,T),\ \mu\in\mathbb{C}.$$

Now we solve the system of equations $(4.1)_R$, $(4.1)_L$, $(4.2)_1$, $(4.2)_2$, for the unknown quantities

$\tilde{g}_1^R(\mu^2,T),\ \tilde{g}_0^R(\mu^2,T),\ \tilde{g}_1^L(\mu^2,T),\ \tilde{g}_0^L(\mu^2,T)$, with $\mu\in\Delta^+$ (fig 7).

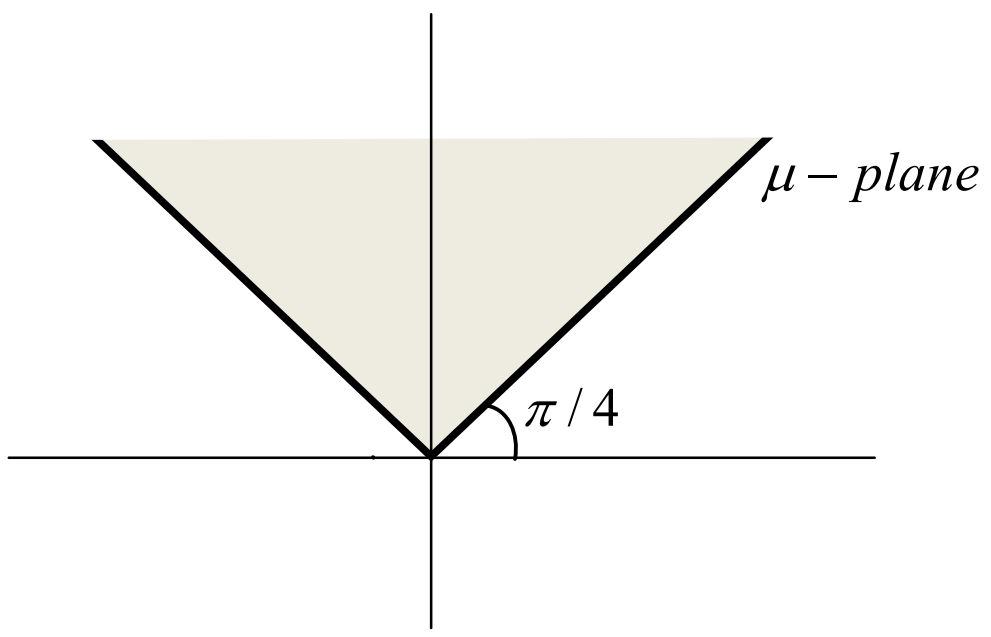

**Fig 7** *The set* $\Delta^+$. *The system (4.1) & (4.2) is solved for* $\mu\in\Delta^+$.

The determinant of the coefficients of the unknowns is

$$\mathfrak{D}(\mu)=\det\begin{bmatrix}1 & -i\varphi_R(\mu) & 0 & 0\\ 0 & 0 & 1 & -i\varphi_L(-\mu)\\ \alpha_{11} & \alpha_{12} & \alpha_{13} & \alpha_{14}\\ \alpha_{21} & \alpha_{22} & \alpha_{23} & \alpha_{24}\end{bmatrix}.$$

By Cramer's rule, we obtain, for $\mu\in\Delta^{+}$,

$$\tilde{g}_1^R(\mu^2,T)=\frac{\mathfrak{D}_1^{R,*}(\mu,T)}{\mathfrak{D}(\mu)},\ \tilde{g}_0^R(\mu^2,T)=\frac{\mathfrak{D}_0^{R,*}(\mu,T)}{\mathfrak{D}(\mu)},\ \tilde{g}_1^L(\mu^2,T)=\frac{\mathfrak{D}_0^{L,*}(\mu,T)}{\mathfrak{D}(\mu)},\ \tilde{g}_0^L(\mu^2,T)=\frac{\mathfrak{D}_0^{L,*}(\mu,T)}{\mathfrak{D}(\mu)},$$

where

$$\mathfrak{D}_1^{R,*}(\mu,T)=\det\begin{bmatrix}J_R^*(\mu,T)+J_R(\mu,T) & -i\varphi_R(\mu) & 0 & 0\\ J_L^*(\mu,T)+J_L(\mu,T) & 0 & 1 & -i\varphi_L(-\mu)\\ \tilde{I}_1(\mu^2,T) & \alpha_{12} & \alpha_{13} & \alpha_{14}\\ \tilde{I}_2(\mu^2,T) & \alpha_{22} & \alpha_{23} & \alpha_{24}\end{bmatrix}$$

$$\mathfrak{D}_0^{R,*}(\mu,T)=\det\begin{bmatrix}1 & J_R^*(\mu,T)+J_R(\mu,T) & 0 & 0\\ 0 & J_L^*(\mu,T)+J_L(\mu,T) & 1 & -i\varphi_L(-\mu)\\ \alpha_{11} & \tilde{I}_1(\mu^2,T) & \alpha_{13} & \alpha_{14}\\ \alpha_{21} & \tilde{I}_2(\mu^2,T) & \alpha_{23} & \alpha_{24}\end{bmatrix}$$

$$\mathfrak{D}_1^{L,*}(\mu,T)=\det\begin{bmatrix}1 & -i\varphi_R(\mu) & J_R^*(\mu,T)+J_R(\mu,T) & 0\\ 0 & 0 & J_L^*(\mu,T)+J_L(\mu,T) & -i\varphi_L(-\mu)\\ \alpha_{11} & \alpha_{12} & \tilde{I}_1(\mu^2,T) & \alpha_{14}\\ \alpha_{21} & \alpha_{22} & \tilde{I}_2(\mu^2,T) & \alpha_{24}\end{bmatrix}$$

$$\mathfrak{D}_0^{L,*}(\mu,T)=\det\begin{bmatrix}1 & -i\varphi_R(\mu) & 0 & J_R^*(\mu,T)+J_R(\mu,T)\\ 0 & 0 & 1 & J_L^*(\mu,T)+J_L(\mu,T)\\ \alpha_{11} & \alpha_{12} & \alpha_{13} & \tilde{I}_1(\mu^2,T)\\ \alpha_{21} & \alpha_{22} & \alpha_{23} & \tilde{I}_2(\mu^2,T)\end{bmatrix}.$$

Let us define

$$G_1^R(\mu,T)=\frac{\mathfrak{D}_1^{R}(\mu,T)}{\mathfrak{D}(\mu)},\ G_0^R(\mu,T)=\frac{\mathfrak{D}_0^{R}(\mu,T)}{\mathfrak{D}(\mu)},\ G_1^L(\mu,T)=\frac{\mathfrak{D}_1^{L}(\mu,T)}{\mathfrak{D}(\mu)},\ G_0^L(\mu,T)=\frac{\mathfrak{D}_0^{L}(\mu,T)}{\mathfrak{D}(\mu)},$$

where $\mathfrak{D}_1^R(\mu,T)$, $\mathfrak{D}_0^R(\mu,T)$, $\mathfrak{D}_1^L(\mu,T)$, $\mathfrak{D}_0^L(\mu,T)$ are the deteninants which result from the corresponding deteninants $\mathfrak{D}_1^{R,*}(\mu,T)$, $\mathfrak{D}_0^{R,*}(\mu,T)$, $\mathfrak{D}_1^{L,*}(\mu,T)$, $\mathfrak{D}_0^{L,*}(\mu,T)$, by eliminating the terms $J_R^*(\mu,T)$ and $J_L^*(\mu,T)$, i.e., by replacing them by zero. Thus,

$$\mathfrak{D}_1^{R,*}(\mu,T)-\mathfrak{D}_1^{R}(\mu,T)=\mathfrak{U}_1^R(\mu,T),\ \mathfrak{D}_0^{R,*}(\mu,T)-\mathfrak{D}_0^{R}(\mu,T)=\mathfrak{U}_0^R(\mu,T),$$

$$\mathfrak{D}_1^{L,*}(\mu,T)-\mathfrak{D}_1^{L}(\mu,T)=\mathfrak{U}_1^L(\mu,T),\ \mathfrak{D}_0^{L,*}(\mu,T)-\mathfrak{D}_0^{L}(\mu,T)=\mathfrak{U}_0^L(\mu,T),$$

where

$$\mathfrak{U}_1^R(\mu,T)=\det\begin{bmatrix} J_R^*(\mu,T) & -i\varphi_R(\mu) & 0 & 0 \\ J_L^*(\mu,T) & 0 & 1 & -i\varphi_L(-\mu) \\ 0 & \alpha_{12} & \alpha_{13} & \alpha_{14} \\ 0 & \alpha_{22} & \alpha_{23} & \alpha_{24} \end{bmatrix},$$

$$\mathfrak{U}_0^R(\mu,T)=\det\begin{bmatrix} 1 & J_R^*(\mu,T) & 0 & 0 \\ 0 & J_L^*(\mu,T) & 1 & -i\varphi_L(-\mu) \\ \alpha_{11} & 0 & \alpha_{13} & \alpha_{14} \\ \alpha_{21} & 0 & \alpha_{23} & \alpha_{24} \end{bmatrix},$$

$$\mathfrak{U}_1^L(\mu,T)=\det\begin{bmatrix} 1 & -i\varphi_R(\mu) & J_R^*(\mu,T) & 0 \\ 0 & 0 & J_L^*(\mu,T) & -i\varphi_L(-\mu) \\ \alpha_{11} & \alpha_{12} & 0 & \alpha_{14} \\ \alpha_{21} & \alpha_{22} & 0 & \alpha_{24} \end{bmatrix},$$

$$\mathfrak{U}_0^L(\mu,T)=\det\begin{bmatrix} 1 & -i\varphi_R(\mu) & 0 & J_R^*(\mu,T) \\ 0 & 0 & 1 & J_L^*(\mu,T) \\ \alpha_{11} & \alpha_{12} & \alpha_{13} & 0 \\ \alpha_{21} & \alpha_{22} & \alpha_{23} & 0 \end{bmatrix}.$$

A computation shows that

$$\mathfrak{D}(\mu)=\mathrm{A}_{13}\varphi_R(\mu)\varphi_L(-\mu)-i\mathrm{A}_{14}\varphi_R(\mu)-i\mathrm{A}_{23}\varphi_L(-\mu)-\mathrm{A}_{24},$$

$$\mathfrak{U}_1^R(\mu,T)=i\varphi_R(\mu)\mathrm{A}_{34}J_L^*(\mu,T)-i\varphi_L(-\mu)\mathrm{A}_{23}J_R^*(\mu,T)-\mathrm{A}_{24}J_R^*(\mu,T),$$

$$\mathfrak{U}_0^R(\mu,T)=i\varphi_L(-\mu)\mathrm{A}_{13}J_R^*(\mu,T)+\mathrm{A}_{14}J_R^*(\mu,T)+\mathrm{A}_{34}J_L^*(\mu,T),$$

$$\mathfrak{U}_1^L(\mu,T)=-i\varphi_R(\mu)\mathrm{A}_{14}J_L^*(\mu,T)-i\varphi_L(-\mu)\mathrm{A}_{12}J_R^*(\mu,T)-\mathrm{A}_{24}J_L^*(\mu,T),$$

$$\mathfrak{U}_0^L(\mu,T)=i\varphi_R(\mu)\mathrm{A}_{13}J_L^*(\mu,T)-\mathrm{A}_{12}J_R^*(\mu,T)+\mathrm{A}_{23}J_L^*(\mu,T),$$

where

$$\mathrm{A}_{13}:=\det\begin{bmatrix}\alpha_{11} & \alpha_{13}\\ \alpha_{21} & \alpha_{23}\end{bmatrix},\ \mathrm{A}_{14}:=\det\begin{bmatrix}\alpha_{11} & \alpha_{14}\\ \alpha_{21} & \alpha_{24}\end{bmatrix},\ \mathrm{A}_{23}:=\det\begin{bmatrix}\alpha_{12} & \alpha_{13}\\ \alpha_{22} & \alpha_{23}\end{bmatrix},$$

$$\mathrm{A}_{24}:=\det\begin{bmatrix}\alpha_{12} & \alpha_{14}\\ \alpha_{22} & \alpha_{24}\end{bmatrix},\ \mathrm{A}_{12}:=\det\begin{bmatrix}\alpha_{11} & \alpha_{12}\\ \alpha_{21} & \alpha_{22}\end{bmatrix}\ \textit{and}\ \mathrm{A}_{34}:=\det\begin{bmatrix}\alpha_{13} & \alpha_{14}\\ \alpha_{23} & \alpha_{24}\end{bmatrix}.$$

***Assumption which guarantees that*** $\mathfrak{D}(\mu)\neq 0$ ***for large*** $|\mu|$**.** We claim that there exists a non-zero polynomial $P(z)$ such that if $\mathfrak{D}(\mu)=0$ for some $\mu\in\mathbb{C}$, then $P(\mu)=0$. Indeed, if $\mathfrak{D}(\mu)=0$ then

$$\mathrm{A}_{13}\varphi_R(\mu)\varphi_L(-\mu)-i\mathrm{A}_{14}\varphi_R(\mu)-i\mathrm{A}_{23}\varphi_L(-\mu)-\mathrm{A}_{24}=0$$

$$\Rightarrow\ \mathrm{A}_{13}[2\sigma_R^2\varphi_R(\mu)][2\sigma_L^2\varphi_L(-\mu)]-2i\mathrm{A}_{14}\sigma_L^2[2\sigma_R^2\varphi_R(\mu)]-2i\mathrm{A}_{23}\sigma_R^2[2\sigma_L^2\varphi_L(-\mu)]-4\sigma_R^2\sigma_L^2\mathrm{A}_{24}=0$$

$$\Rightarrow\ \mathrm{A}_{13}[i\gamma_R+2\sigma_R\sqrt[R]{\mu^2-\varepsilon_R^2}][i\gamma_L-2\sigma_L\sqrt[L]{\mu^2-\varepsilon_L^2}]$$

$$-2i\mathrm{A}_{14}\sigma_L^2[i\gamma_R+2\sigma_R\sqrt[R]{\mu^2-\varepsilon_R^2}]-2i\mathrm{A}_{23}\sigma_R^2[i\gamma_L-2\sigma_L\sqrt[L]{\mu^2-\varepsilon_L^2}]-4\sigma_R^2\sigma_L^2\mathrm{A}_{24}=0$$

$$\Rightarrow \mathcal{B}_1 \sqrt[R]{\mu^2-\varepsilon_R^2}\sqrt[L]{\mu^2-\varepsilon_L^2}+\mathcal{B}_2\sqrt[L]{\mu^2-\varepsilon_L^2}=\mathcal{B}_3\sqrt[R]{\mu^2-\varepsilon_R^2}+\mathcal{B}_4,$$

where

$$\mathcal{B}_1=-4\mathrm{A}_{13}\sigma_R\sigma_L,$$

$$\mathcal{B}_2=-2i\mathrm{A}_{13}\sigma_L\gamma_R+4i\mathrm{A}_{23}\sigma_L\sigma_R^2,$$

$$\mathcal{B}_3=-2i\mathrm{A}_{13}\sigma_R\gamma_L+4i\mathrm{A}_{14}\sigma_R\sigma_L^2,$$

$$\mathcal{B}_4=\mathrm{A}_{13}\gamma_R\gamma_L-2\mathrm{A}_{14}\gamma_R\sigma_L^2-2\mathrm{A}_{23}\gamma_L\sigma_R^2+4\mathrm{A}_{24}\sigma_R^2\sigma_L^2.$$

Thus $\mathfrak{D}(\mu)=0 \;\Rightarrow\; \mathcal{B}_1^2(\mu^2-\varepsilon_R^2)(\mu^2-\varepsilon_L^2)+\mathcal{B}_2^2(\mu^2-\varepsilon_L^2)+2\mathcal{B}_1\mathcal{B}_2(\mu^2-\varepsilon_L^2)\sqrt[R]{\mu^2-\varepsilon_R^2}$

$$=\mathcal{B}_3^2(\mu^2-\varepsilon_R^2)+2\mathcal{B}_3\mathcal{B}_4\sqrt[R]{\mu^2-\varepsilon_R^2}+\mathcal{B}_4^2 \;\Rightarrow\; P(\mu)=0,$$

where we have set

$$P(\mu):=[\mathcal{B}_1^2(\mu^2-\varepsilon_R^2)(\mu^2-\varepsilon_L^2)+\mathcal{B}_2^2(\mu^2-\varepsilon_L^2)-\mathcal{B}_3^2(\mu^2-\varepsilon_R^2)-\mathcal{B}_4^2]^2-4(\mu^2-\varepsilon_R^2)[\mathcal{B}_1\mathcal{B}_2(\mu^2-\varepsilon_L^2)+\mathcal{B}_3\mathcal{B}_4]^2.$$

*Remarks* (1) We note that if $\mathrm{A}_{13}\neq 0$ then the degree of the polynomial $P(\mu)$ is $8$, and, in particular, $P(\mu)\not\equiv 0$.

(2) Similarly if $\mathrm{A}_{13}=0$ and $\mathrm{A}_{23}\sigma_R\neq\pm\mathrm{A}_{14}\sigma_L$ then the degree of the polynomial $P(\mu)$ is $4$.

(3) In the case $\mathrm{A}_{12}=\mathrm{A}_{34}=0$ and $\mathrm{A}_{24}\neq 0$, the interface problem [$(1.1)_R$, $(1.1)_L$, (1.2)] is "separable" in the sense that it reduces to two independant half-line problems for the heat equation.

**Theorem** *Suppose that* $P(\mu)\not\equiv 0$. *Then the function* $\mathfrak{D}(\mu)$ *has finitely many zeros and, in particular, there exists* $\rho>0$ *such that* $\mathfrak{D}(\mu)\neq 0$ *for every* $\mu\in\mathbb{C}$, *with* $|\mu|\geq\rho$.
*For such a* $\rho$, *the UTM solution*

$$\{[u^R(x,t),\, x>0,\, t>0],\, [u^L(x,t),\, x<0,\, t>0]\},$$

*of the interface problem {* $(1.1)_R$, $(1.1)_L$ *,(1.2)}, given by (2.8) and (3.8), can be written in the following way:*

$$(4.3)\quad u^R(x,t)=\frac{1}{2\pi}\int_{-\infty}^{\infty}e^{i\lambda x-\omega_R(\lambda)t}\hat{u}_0^R(\lambda)d\lambda$$
$$-\frac{1}{2\pi}\int_{\mu\in\partial \mathrm{E}_\rho^+}e^{i\varphi_R(\mu)x-\mu^2 t}\{\sigma_R^2G_1^R(\mu,T)+[\sigma_R^2 i\varphi_R(\mu)+\gamma_R]G_0^R(\mu,T)\}\frac{d\varphi_R(\mu)}{d\mu}d\mu$$
$$+\frac{1}{2\pi}\int_{-\infty}^{\infty}e^{i\lambda x-\omega_R(\lambda)t}\tilde{\hat{f}}_R(\lambda,\omega_R,t)d\lambda,\ \text{for } x>0,\ 0<t<T,$$

$$(4.4)\quad u^L(x,t)=\frac{1}{2\pi}\int_{-\infty}^{\infty}e^{i\lambda x-\omega_L(\lambda)t}\hat{u}_0^L(\lambda)d\lambda$$
$$+\frac{1}{2\pi}\int_{\mu\in\partial \mathrm{E}_\rho^+}e^{i\varphi_L(-\mu)x-\mu^2 t}\{\sigma_L^2G_1^L(\mu,T)+[\sigma_L^2 i\varphi_L(-\mu)+\gamma_L]G_0^L(\mu,T)\}\left\{-\frac{d[\varphi_L(-\mu)]}{d\mu}\right\}d\mu$$
$$+\frac{1}{2\pi}\int_{-\infty}^{\infty}e^{i\lambda x-\omega_L(\lambda)t}\tilde{\hat{f}}_L(\lambda,\omega_L,t)d\lambda,\ \text{for } x<0,\ 0<t<T,$$

*where* (fig 8)

$$\mathrm{E}_\rho^+:=\{\mu\in\mathbb{C}:\ \mathrm{Im}\,\mu\geq 0,\ \mathrm{Re}(\mu^2)\geq 0\ \textit{and}\ |\mu|\geq\rho\}.$$

*Also we have*

$$(4.5)\, u^R(x,t)=\frac{1}{2\pi}\int_{-\infty}^{\infty}e^{i\lambda x-\omega_R(\lambda)t}\hat{u}_0^R(\lambda)d\lambda$$

$$-\frac{1}{2\pi}\int_{\lambda\in\partial\mathcal{Z}^{+}_{\rho_0}} e^{i\lambda x-\omega_R(\lambda)t}\{\sigma_R^2 G_1^R(\psi_R(\lambda),T)+[\sigma_R^2 i\lambda+\gamma_R]G_0^R(\psi_R(\lambda),T)\}d\lambda$$

$$+\frac{1}{2\pi}\int_{-\infty}^{\infty} e^{i\lambda x-\omega_R(\lambda)t}\tilde{\hat{f}}_R(\lambda,\omega_R,t)d\lambda\,, \textit{for } x>0\,,\ 0<t<T\,,$$

*and*

(4.6) $$u^L(x,t)=\frac{1}{2\pi}\int_{-\infty}^{\infty} e^{i\lambda x-\omega_L(\lambda)t}\hat{u}_0^L(\lambda)d\lambda$$

$$+\frac{1}{2\pi}\int_{\lambda\in\partial\mathcal{Z}^{-}_{\rho_0}} e^{i\lambda x-\omega_L(\lambda)t}\{\sigma_L^2 G_1^L(-\psi_L(\lambda),T)+[\sigma_L^2 i\lambda+\gamma_L]G_0^L(-\psi_R(\lambda),T)\}d\lambda$$

$$+\frac{1}{2\pi}\int_{-\infty}^{\infty} e^{i\lambda x-\omega_L(\lambda)t}\tilde{\hat{f}}_L(\lambda,\omega_L,t)d\lambda\,, \text{ for } x<0\,,\ 0<t<T\,,$$

*where* (fig. 9 and fig. 10)

$$\mathcal{Z}^{+}_{\rho_0}=\{\lambda\in\mathbb{C}:\operatorname{Im}\lambda\geq 0,\ \operatorname{Re}[\omega_R(\lambda)]\geq 0 \textit{ and } |\lambda|\geq\rho_0\}$$

*and*

$$\mathcal{Z}^{-}_{\rho_0}=\{\lambda\in\mathbb{C}:\operatorname{Im}\lambda\leq 0,\ \operatorname{Re}[\omega_L(\lambda)]\geq 0 \textit{ and } |\lambda|\geq\rho_0\}\,,$$

*provided that* $\rho_0>0$ *is sufficiently large.*

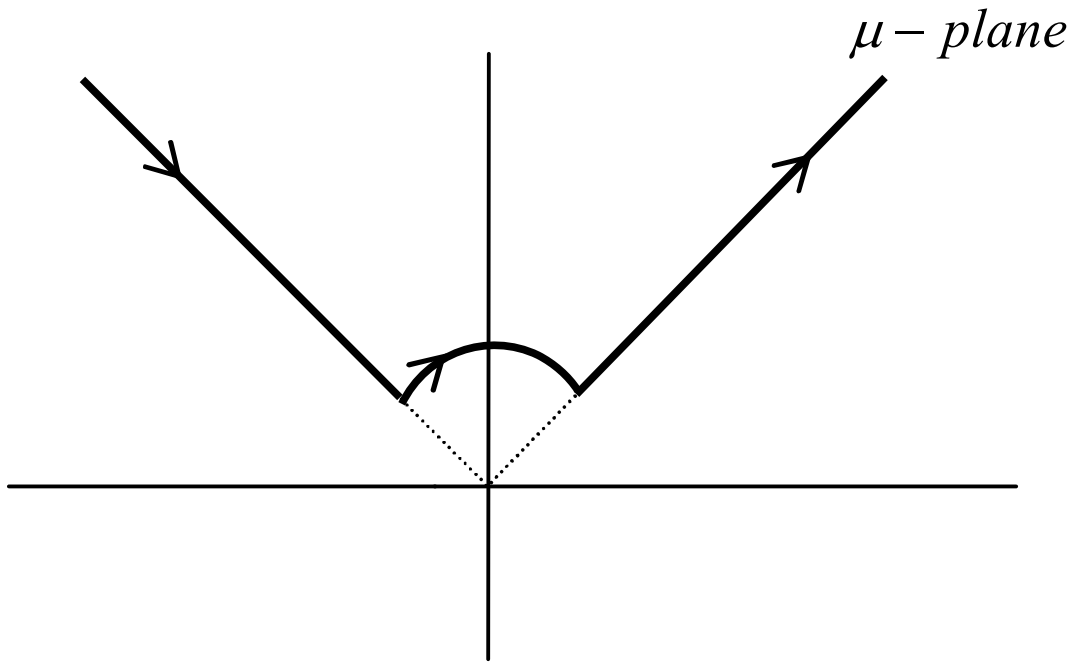

**Fig 8** *The contour* $\partial\mathrm{E}^{+}_{\rho}$

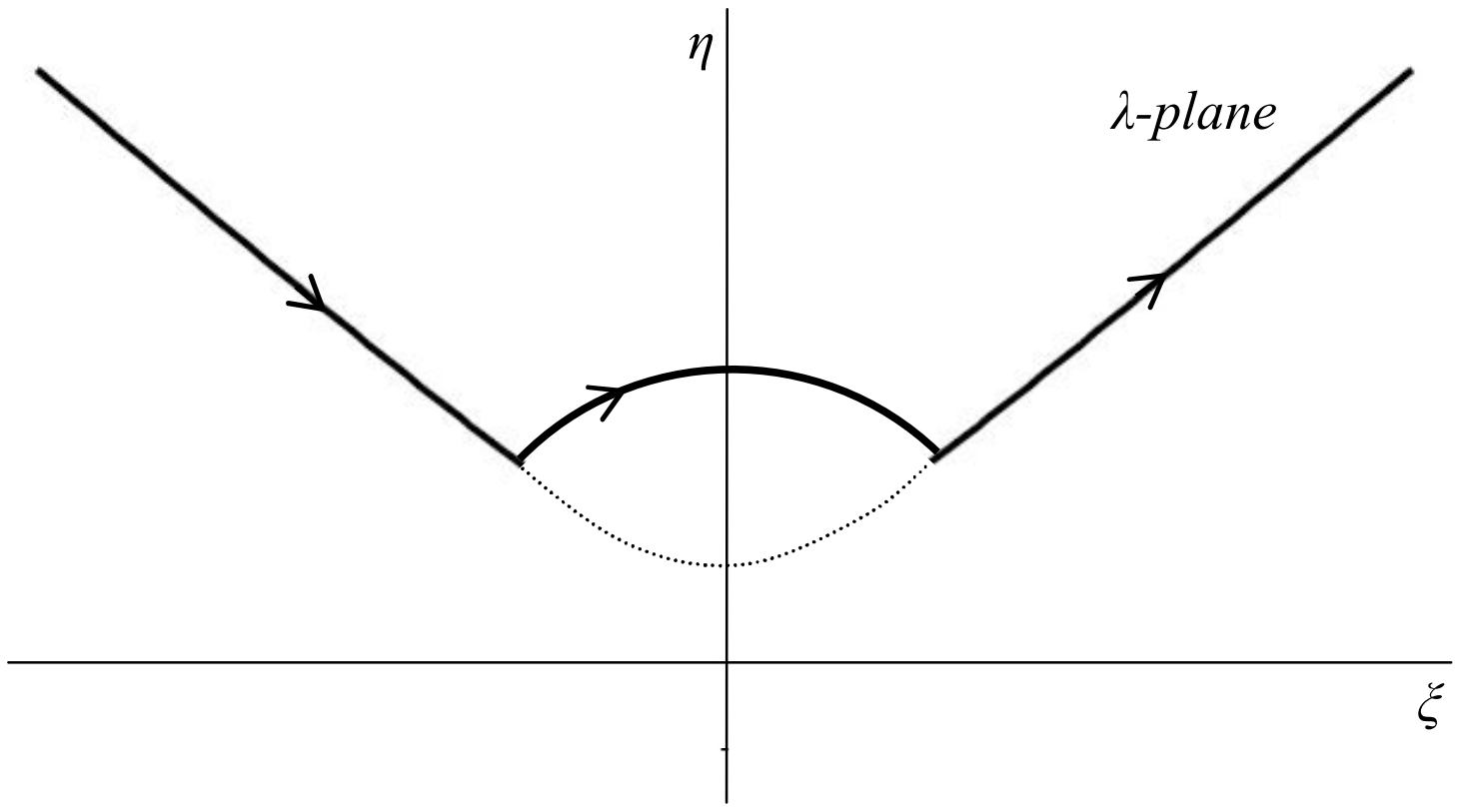

**Fig 9** *The contour* $\partial\mathcal{Z}^{+}_{\rho_0}$

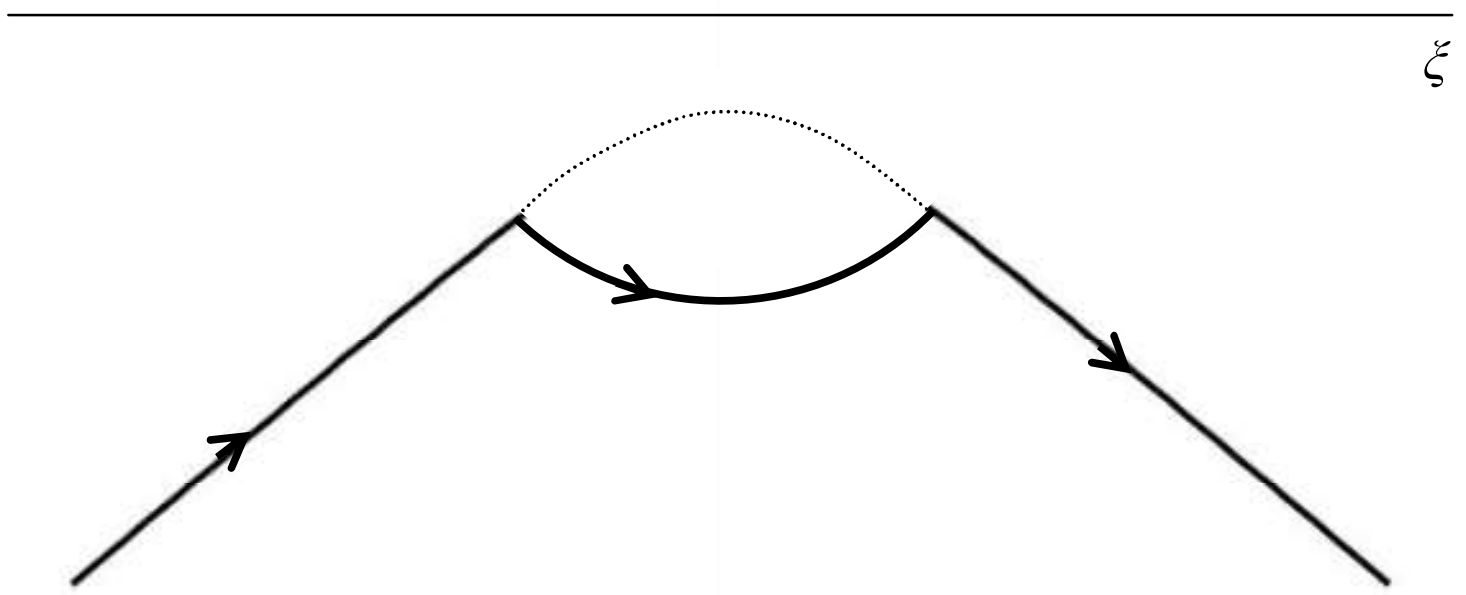

**Fig 10** *The contour* $\partial\mathcal{Z}^{-}_{\rho_0}$

**Proof** Substituting $\tilde{g}_1^R(\mu^2,T),\ \tilde{g}_0^R(\mu^2,T),\ \tilde{g}_1^L(\mu^2,T),\ \tilde{g}_0^L(\mu^2,T)$ in (2.8) and (3.8), and deforming the contour $\Gamma^+$ to the contour $\partial\mathrm{E}^+_\rho$, in order to avoid the zeros of the function $\det\mathfrak{D}(\mu)$, we obtain the following formulas:

$$(4.7)\quad u^R(x,t)=\frac{1}{2\pi}\int_{-\infty}^{\infty}e^{i\lambda x-\omega_R(\lambda)t}\hat{u}_0^R(\lambda)d\lambda$$
$$-\frac{1}{2\pi}\int_{\mu\in\partial\mathrm{E}^+_\rho}e^{i\varphi_R(\mu)x-\mu^2 t}\left\{\sigma_R^2\frac{\mathfrak{D}_1^{R,*}(\mu,T)}{\mathfrak{D}(\mu)}+[\sigma_R^2 i\varphi_R(\mu)+\gamma_R]\frac{\mathfrak{D}_0^{R,*}(\mu,T)}{\mathfrak{D}(\mu)}\right\}\frac{d\varphi_R(\mu)}{d\mu}d\mu$$
$$+\frac{1}{2\pi}\int_{-\infty}^{\infty}e^{i\lambda x-\omega_R(\lambda)t}\tilde{\tilde{f}}_R(\lambda,\omega_R,t)d\lambda,\ \text{for}\ x>0,\ 0<t<T,$$

and

$$(4.8)\quad u^L(x,t)=\frac{1}{2\pi}\int_{-\infty}^{\infty}e^{i\lambda x-\omega_L(\lambda)t}\hat{u}_0^L(\lambda)d\lambda$$
$$+\frac{1}{2\pi}\int_{\mu\in\partial\mathrm{E}^+_\rho}e^{i\varphi_L(-\mu)x-\mu^2 t}\left\{\sigma_L^2\frac{\mathfrak{D}_1^{L,*}(\mu,T)}{\mathfrak{D}(\mu)}+[\sigma_L^2 i\varphi_L(-\mu)+\gamma_L]\frac{\mathfrak{D}_0^{L,*}(\mu,T)}{\mathfrak{D}(\mu)}\right\}\left\{-\frac{d[\varphi_L(-\mu)]}{d\mu}\right\}d\mu$$
$$+\frac{1}{2\pi}\int_{-\infty}^{\infty}e^{i\lambda x-\omega_L(\lambda)t}\tilde{\tilde{f}}_L(\lambda,\omega_L,t)d\lambda,\ \text{for}\ x<0,\ 0<t<T.$$

Since the coefficients of the quantities $J_R^*(\mu,T)$, $J_L^*(\mu,T)$, in the determinants $\mathfrak{U}_1^R(\mu,T)$, $\mathfrak{U}_0^R(\mu,T)$, $\mathfrak{U}_1^L(\mu,T)$, $\mathfrak{U}_0^L(\mu,T)$, are $\mathrm{O}(\mu)$, as $\mu\to\infty$ with $\operatorname{Im}\mu\geq 0$, it follows that

$$\mathfrak{U}_1^R(\mu,T)=\mathrm{O}(e^{\mu^2 T}),\ \mathfrak{U}_0^R(\mu,T)=\mathrm{O}(e^{\mu^2 T}),\ \mathfrak{U}_1^L(\mu,T)=\mathrm{O}(e^{\mu^2 T}),\ \mathfrak{U}_0^L(\mu,T)=\mathrm{O}(e^{\mu^2 T}).$$

Also, $\mathfrak{D}(\mu)=\mathrm{O}(\mu^2)$ and, therefore,

$$(4.9)\quad \tilde{g}_0^R(\mu^2,T)-G_0^R(\mu,T)=\mathrm{O}(e^{\mu^2 T}),\ \tilde{g}_1^R(\mu^2,T)-G_1^R(\mu,T)=\mathrm{O}(e^{\mu^2 T})$$

and

$$(4.10)\quad \tilde{g}_0^L(\mu^2,T)-G_0^L(\mu,T)=\mathrm{O}(e^{\mu^2 T}),\ \tilde{g}_1^L(\mu^2,T)-G_1^L(\mu,T)=\mathrm{O}(e^{\mu^2 T}).$$

Also, keeping in mind that

$$\frac{d\varphi_R(\mu)}{d\mu}=\mathrm{O}(1),\ \frac{d\varphi_L(-\mu)}{d\mu}=\mathrm{O}(1),$$

we obtain

(4.11) $$\int_{\mu\in\partial\mathrm{E}_\rho^+} e^{i\varphi_R(\mu)x-\mu^2 t}\left\{\sigma_R^2\frac{\mathfrak{U}_1^R(\mu,T)}{\mathfrak{D}(\mu)}+[\sigma_R^2 i\varphi_R(\mu)+\gamma_R]\frac{\mathfrak{U}_0^R(\mu,T)}{\mathfrak{D}(\mu)}\right\}\frac{d\varphi_R(\mu)}{d\mu}d\mu=0, \text{ for } x>0,$$

and

(4.12) $$\int_{\mu\in\partial\mathrm{E}_\rho^+} e^{i\varphi_L(-\mu)x-\mu^2 t}\left\{\sigma_L^2\frac{\mathfrak{U}_1^L(\mu,T)}{\mathfrak{D}(\mu)}+[\sigma_L^2 i\varphi_L(-\mu)+\gamma_L]\frac{\mathfrak{U}_0^L(\mu,T)}{\mathfrak{D}(\mu)}\right\}\left\{-\frac{d[\varphi_L(-\mu)]}{d\mu}\right\}d\mu=0, \text{ for } x<0.$$

Indeed, the above equations follow from (4.9), (4.10), Cauchy's theorem and the following form of Jordan's lemma:

$$\lim_{\mathrm{B}\to+\infty}\int_{\{\mu:\,|\mu|=\mathrm{B}\}\cap\mathrm{E}_\rho^+} e^{i\varphi_R(\mu)x+\mu^2(T-t)}\mu^N d\mu=0, \text{ for } x>0 \text{ and } N\in\mathbb{N}.$$

*Completion of the proof of the theorem* In view of equations (4.11) and (4.12), formulas (4.3) and (4.4) follow from (4.7) and (4.8). Finally, changing the variable form $\mu$ to $\lambda$, we obtain (4.5) and (4.6).

*Remark* Formulas (4.5) and (4.6) hold also if we set $T$ in place of $t$, i.e., if we replace the quantities

$$G_1^R(\psi_R(\lambda),T),\ G_0^R(\psi_R(\lambda),T),\ G_1^L(-\psi_L(\lambda),T),\ G_0^L(-\psi_L(\lambda),T)$$

by

$$G_1^R(\psi_R(\lambda),t),\ G_0^R(\psi_R(\lambda),t),\ G_1^L(-\psi_L(\lambda),t),\ G_0^L(-\psi_L(\lambda),t).$$

In conclusion, the principal aim of this paper was to derive, for the first time, explicit solution formulae for a broad class of one-dimensional (representing a symmetry in higher-dimensional configurations) interface problems associated with the lRAD equation (with discontinuous coefficients and fairly arbitrary interface conditions). Our approach combined the Fokas' unified transform methodology with complex-analytic techniques in order to obtain integral representations of the solution directly in terms of the prescribed initial, interface and forcing data. In particular, the analysis yielded exact characterisations of the solution operators which are expected to clarify the role played by the interface conditions in determining the qualitative behavior of the associated evolution. The presence of the lower-order terms in the lRAD operator (as opposed to the traditional heat equation), together with the general form of interface conditions that we addressed, introduced structural features and analytical complications absent from the existing literature.

The general analytical framework developed herein opens the way to a systematic investigation of a remarkably broad class of heterogeneous transport phenomena, extending well beyond the classical diffusion and perfect-contact interface models. In forthcoming work, we will demonstrate the versatility and analytical power of the present solution representations through physically motivated examples involving imperfect transmission, interfacial permeability and reactivity, and heterogeneous advection–diffusion–reaction dynamics, together with comparisons against established transform-based approaches. Particular emphasis will be placed on the explicit characterization of interfacial fluxes, transmission and reflection mechanisms, and asymptotic regimes, thereby elucidating the broader applicability and advantages of the present methodology.

## Acknowledgment

This research work was implemented in the framework of H.F.R.I call "Basic research Financing (Horizontal support of all Sciences)" under the National Recovery and Resilience Plan "Greece 2.0" funded by the European Union – NextGenerationEU (H.F.R.I. Project Number: 14910), for which AC (supported until 11/2025) and GK express their thankfulness. AC is currently supported by JST [Moonshot R&D][Grant Number JPMJMS25A5].